\PassOptionsToPackage{unicode}{hyperref}
\PassOptionsToPackage{hyphens}{url}
\PassOptionsToPackage{dvipsnames,svgnames,x11names}{xcolor}
\documentclass[
  12pt,letterpaper]{article}

\usepackage{amsmath,amssymb}
\usepackage{iftex}
\ifPDFTeX
  \usepackage[T1]{fontenc}
  \usepackage[utf8]{inputenc}
  \usepackage{textcomp} 
\else 
  \usepackage{unicode-math}
  \defaultfontfeatures{Scale=MatchLowercase}
  \defaultfontfeatures[\rmfamily]{Ligatures=TeX,Scale=1}
\fi
\usepackage{lmodern}
\ifPDFTeX\else  
\fi
\IfFileExists{upquote.sty}{\usepackage{upquote}}{}
\IfFileExists{microtype.sty}{
  \usepackage[]{microtype}
  \UseMicrotypeSet[protrusion]{basicmath} 
}{}
\makeatletter
\@ifundefined{KOMAClassName}{
  \IfFileExists{parskip.sty}{%
    \usepackage{parskip}
  }{
    \setlength{\parindent}{0pt}
    \setlength{\parskip}{6pt plus 2pt minus 1pt}}
}{
  \KOMAoptions{parskip=half}}
\makeatother
\usepackage{xcolor}
\makeatletter
\ifx\paragraph\undefined\else
  \let\oldparagraph\paragraph
  \renewcommand{\paragraph}{
    \@ifstar
      \xxxParagraphStar
      \xxxParagraphNoStar
  }
  \newcommand{\xxxParagraphStar}[1]{\oldparagraph*{#1}\mbox{}}
  \newcommand{\xxxParagraphNoStar}[1]{\oldparagraph{#1}\mbox{}}
\fi
\ifx\subparagraph\undefined\else
  \let\oldsubparagraph\subparagraph
  \renewcommand{\subparagraph}{
    \@ifstar
      \xxxSubParagraphStar
      \xxxSubParagraphNoStar
  }
  \newcommand{\xxxSubParagraphStar}[1]{\oldsubparagraph*{#1}\mbox{}}
  \newcommand{\xxxSubParagraphNoStar}[1]{\oldsubparagraph{#1}\mbox{}}
\fi
\makeatother

\usepackage{longtable,booktabs,array}
\usepackage{calc} 
\usepackage{etoolbox}
\makeatletter
\patchcmd\longtable{\par}{\if@noskipsec\mbox{}\fi\par}{}{}
\makeatother
\IfFileExists{footnotehyper.sty}{\usepackage{footnotehyper}}{\usepackage{footnote}}
\makesavenoteenv{longtable}
\usepackage{graphicx}
\makeatletter
\def\maxwidth{\ifdim\Gin@nat@width>\linewidth\linewidth\else\Gin@nat@width\fi}
\def\maxheight{\ifdim\Gin@nat@height>\textheight\textheight\else\Gin@nat@height\fi}
\makeatother
\setkeys{Gin}{width=\maxwidth,height=\maxheight,keepaspectratio}
\makeatletter
\def\fps@figure{htbp}
\makeatother

\makeatletter
\@ifpackageloaded{caption}{}{\usepackage{caption}}
\AtBeginDocument{%
\ifdefined\contentsname
  \renewcommand*\contentsname{Table of contents}
\else
  \newcommand\contentsname{Table of contents}
\fi
\ifdefined\listfigurename
  \renewcommand*\listfigurename{List of Figures}
\else
  \newcommand\listfigurename{List of Figures}
\fi
\ifdefined\listtablename
  \renewcommand*\listtablename{List of Tables}
\else
  \newcommand\listtablename{List of Tables}
\fi
\ifdefined\figurename
  \renewcommand*\figurename{Figure}
\else
  \newcommand\figurename{Figure}
\fi
\ifdefined\tablename
  \renewcommand*\tablename{Table}
\else
  \newcommand\tablename{Table}
\fi
}
\@ifpackageloaded{float}{}{\usepackage{float}}
\floatstyle{ruled}
\@ifundefined{c@chapter}{\newfloat{codelisting}{h}{lop}}{\newfloat{codelisting}{h}{lop}[chapter]}
\floatname{codelisting}{Listing}

\makeatother
\makeatletter
\@ifpackageloaded{caption}{}{\usepackage{caption}}
\@ifpackageloaded{subcaption}{}{\usepackage{subcaption}}
\makeatother

\ifLuaTeX
  \usepackage{selnolig}  
\fi
\usepackage[]{natbib}
\usepackage{hyperref}
\usepackage{bookmark}

\IfFileExists{xurl.sty}{\usepackage{xurl}}{} 
\hypersetup{
  pdftitle={Smoothed Least-Absolute Conditional Multidimensional Scaling via Iteratively Reweighted Majorization},
  pdfauthor={Xiao Ling; Anh T. Bui},
  pdfkeywords={Robust estimation, iteratively reweighted least squares, pairwise contamination, nonconvex optimization, distance geometry},
  colorlinks=true,
  linkcolor={blue},
  filecolor={Maroon},
  citecolor={Blue},
  urlcolor={Blue},
  pdfcreator={LaTeX}}

\newcommand{\anon}{1}

\usepackage{amsthm,mathtools,needspace}
\forcsvlist{\listadd\manuscriptenvs}{theorem,proposition,lemma,corollary,assumption}
\newcommand{\protecttheorem}[1]{\AtBeginEnvironment{#1}{\Needspace{4\baselineskip}}}
\forlistloop{\protecttheorem}{\manuscriptenvs}
\newtheorem{theorem}{Theorem}[section]
\newtheorem{proposition}[theorem]{Proposition}

\newtheorem{assumption}[theorem]{Assumption}
\theoremstyle{definition}

\theoremstyle{remark}

\usepackage{algorithmic,algorithm,tabularx,placeins}
\usepackage{tikz,pgfplots}
\usepgfplotslibrary{groupplots}
\pgfplotsset{compat=1.18}

\newcommand{\R}{\mathbb{R}}

\newcommand{\tr}{\operatorname{tr}}

\begin{document}
\def\spacingset#1{\renewcommand{\baselinestretch}{#1}\small\normalsize}
\spacingset{1}
\if1\anon
{
\title{\bfseries Smoothed Least-Absolute Conditional Multidimensional Scaling via Iteratively Reweighted Majorization}
\author{Xiao Ling\thanks{Funding and acknowledgements to be supplied.}\\
Department of Mathematics, Auburn University at Montgomery\\
7400 East Dr, Montgomery, AL 36117, USA\\
\texttt{xling@aum.edu}}
\date{}\maketitle
}\fi
\if0\anon
{\bigskip\bigskip\bigskip\begin{center}
{\LARGE\bfseries Smoothed Least-Absolute Conditional Multidimensional Scaling via Iteratively Reweighted Majorization}
\end{center}\medskip}\fi
\begin{abstract}
Conditional multidimensional scaling combines latent coordinates with observed covariates to fit pairwise dissimilarities. We develop a smoothed $\ell_1$ formulation using the shifted Charbonnier loss. We prove residualwise sharpness of the quadratic loss majorizer, monotonicity, compact containment, and stationarity of noncollapsed accumulation points. A uniform additive bound quantifies approximation to exact least-absolute conditional stress. The experiments use paired contamination, separately calibrated loss parameters, combined-configuration recovery, and direct checks of the theorem's inequalities. 
\end{abstract}
\vfill\newpage\spacingset{1.8}

\section{Introduction}
Pairwise dissimilarities often provide a more natural description of objects than measured coordinates. Examples include similarity judgments in perception studies and comparisons of stochastic surfaces \citep{jaworska2009review,Bui27052021}. Multidimensional scaling (MDS) constructs a low-dimensional configuration whose distances approximate such observations. When some relevant features are already observed, conditional MDS augments them with latent coordinates and estimates their relative contribution to the fitted geometry \citep{bui2023dimension}. 

Squared conditional stress can be dominated by a small number of corrupted dissimilarities. Least-absolute residual fitting reduces the growth of individual contributions, but its usual quadratic reweighting is singular at an exactly fitted residual \citep{heiser1988multidimensional,holland1977robust}. Fair-loss rcMDS resolves this difficulty through finite weights and conditional SMACOF updates \citep{ling2026robust}. We take that model and majorization framework as the starting point and replace the Fair criterion by the shifted Charbonnier loss
\begin{equation}
\rho_\varepsilon(r)=\sqrt{r^2+\varepsilon^2}-\varepsilon,
\qquad \varepsilon>0.
\label{eq:loss}
\end{equation}
The loss is classical \citep{charbonnier1997deterministic}. Its role here is to provide both finite reweighting and a uniform approximation to absolute residuals.

The main analysis follows the rcMDS sequence of residual majorization, monotonicity, compactness, and stationarity. The required adaptation is substantive: the weights become $w_k/\sqrt{r_k^2+\varepsilon^2}$, the uniform curvature lower bound changes accordingly, and the gradient must be related to the \emph{first} exact inner update when more than one update is accepted. Residual smoothing does not make the distance map differentiable at a collapsed pair with positive dissimilarity. We therefore state the needed separation condition at accumulation points explicitly.

The scalar penalty, iteratively reweighted least squares, and geometric SMACOF majorization are established ideas \citep{de2009multidimensional,lange2000optimization,hunter2004tutorial}. Our contribution is a complete conditional-model analysis that connects them: a sharp scalar majorizer, compact centered sublevels, quantitative first-step decrease, qualified stationary-point convergence, and an exact-$\ell_1$ approximation guarantee. The experiments test the inequalities used in the proof and compare recovery under a single documented protocol. Because the two coordinate blocks need not be separately identifiable, the principal recovery target is the combined fitted configuration or its distances.

Section~\ref{sec:model} defines the model and residual majorization. Section~\ref{sec:algorithm} gives the algorithm and convergence results. Section~\ref{sec:approximation} establishes the least-absolute approximation bound. Section~\ref{sec:experiments} presents the computational study, followed by limitations and conclusions.

\section{Model and Residual Majorization}\label{sec:model}
Let $\Delta=(\delta_{ij})$ be a symmetric nonnegative dissimilarity matrix with zero diagonal, and let $\mathbf W=(w_{ij})$ be symmetric and nonnegative. Define the active pairs
\[
\mathcal P=\{(i,j):1\le i<j\le n,\ w_{ij}>0\}.
\]
Asymmetric observations are symmetrized by averaging their two directed values. Write $k=(i,j)$, $\delta_k=\delta_{ij}$, and $w_k=w_{ij}$. Given observed covariates $\mathbf V\in\R^{n\times q}$, latent coordinates $\mathbf U\in\R^{n\times p}$, and a transformation $\mathbf B\in\R^{q\times q}$, set
\begin{align}
\mathbf X&=(\mathbf U,\mathbf B),\qquad
\mathbf Z(\mathbf X)=[\mathbf U,\mathbf V\mathbf B],\nonumber\\
d_k(\mathbf X)&=\left(\|\mathbf u_i-\mathbf u_j\|_2^2+
\|\mathbf B^\top(\mathbf v_i-\mathbf v_j)\|_2^2\right)^{1/2},\label{eq:distance}\\
r_k(\mathbf X)&=\delta_k-d_k(\mathbf X).
\end{align}
The proposed objective is
\begin{equation}
\sigma_\varepsilon(\mathbf X)=\sum_{k\in\mathcal P}w_k\rho_\varepsilon(r_k(\mathbf X)),
\qquad \min_{\mathbf X\in\Omega_0}\sigma_\varepsilon(\mathbf X),
\label{eq:objective}
\end{equation}
where
\[
\Omega_0=\{(\mathbf U,\mathbf B)\in\R^{n\times p}\times\R^{q\times q}:\mathbf1^\top\mathbf U=0\}.
\]
The transformation is unrestricted; no diagonal or orthogonality constraint is imposed. The case $q=0$ omits the covariate block and gives ordinary robust MDS.

\begin{assumption}[Active graph and covariates]\label{ass:regularity}
The graph $(\{1,\ldots,n\},\mathcal P)$ is connected. The covariates are centered, $\mathbf1^\top\mathbf V=0$, and have full column rank. Conditions involving $\mathbf V$ are omitted when $q=0$.
\end{assumption}
Centering leaves all covariate differences unchanged. Under connectedness, full column rank is equivalent to the active differences $\{\mathbf v_i-\mathbf v_j:k\in\mathcal P\}$ spanning $\R^q$. These conditions are needed for compactness and unique centered quadratic updates, but not for residual majorization or monotonicity.

\subsection{Interpretation and loss properties}
Both coordinate blocks admit orthogonal changes of basis, and the latent block admits translation before centering. More fundamentally, when $p,q\ge1$ and $\mathbf b\ne0$, the pairs
\[
(\mathbf V\mathbf b(\mathbf e_1^{(p)})^\top,0)
\quad\hbox{and}\quad
(0,\mathbf b(\mathbf e_1^{(q)})^\top)
\]
have identical distances. Thus conditioning alone does not uniquely separate latent from covariate contributions. We evaluate the combined geometry and interpret latent displays descriptively.

The score and curvature of \eqref{eq:loss} are
\[
\rho_\varepsilon'(r)=\frac{r}{\sqrt{r^2+\varepsilon^2}},
\qquad
\rho_\varepsilon''(r)=\frac{\varepsilon^2}{(r^2+\varepsilon^2)^{3/2}}.
\]
The loss is nonnegative, even, strictly convex in its scalar argument, and coercive, with $|\rho_\varepsilon'(r)|\le1$. Bounded residual score limits the contribution of an individual residual relative to squared loss; it is not a breakdown-point or influence-function guarantee for the fitted configuration. Moreover, \eqref{eq:objective} remains nonconvex in its coordinates and can be nonsmooth at positive-dissimilarity collisions.

\begin{theorem}[Computational hardness]\label{thm:hardness}
For every fixed $\varepsilon>0$, global minimization of \eqref{eq:objective} is NP-hard in general, already for $p=1$, $q=0$, and unit pair weights.
\end{theorem}
\begin{proof}
Use the one-dimensional embedding reduction of \citet[Section~5.1.1]{cayton2006robust} with $f(x)=|x|$ and $g(x)=\rho_\varepsilon(x)$. Both are even and increasing on $[0,\infty)$. For $0\le y\le x\le1$,
\[
\frac{(x-y)x}{2\sqrt{1+\varepsilon^2}}
\le g(x)-g(y)
=\frac{x^2-y^2}{\sqrt{x^2+\varepsilon^2}+\sqrt{y^2+\varepsilon^2}}
\le x-y.
\]
Also $(x-y)x\le f(x)-f(y)\le x-y$. Hence the reduction's upper and lower bounds hold with $\lambda_U=1$ and $\lambda_L=1/(2\sqrt{1+\varepsilon^2})$. Translation to the centered gauge does not change the objective.
\end{proof}

\Needspace{8\baselineskip}
\subsection{Quadratic majorization and residualwise sharpness}
\begin{theorem}[Residual majorization]\label{thm:majorization}
At a support point $\mathbf Y$, define
\begin{align}
q_k(\mathbf Y)&=\frac{w_k}{\sqrt{r_k(\mathbf Y)^2+\varepsilon^2}},\label{eq:weights}\\
Q(\mathbf X\mid\mathbf Y)&=\sigma_\varepsilon(\mathbf Y)
+\frac12\sum_{k\in\mathcal P}q_k(\mathbf Y)
\{r_k(\mathbf X)^2-r_k(\mathbf Y)^2\}.
\label{eq:outer_majorizer}
\end{align}
Then $\sigma_\varepsilon(\mathbf X)\le Q(\mathbf X\mid\mathbf Y)$ for every $\mathbf X$, with equality at $\mathbf X=\mathbf Y$. Minimizing $Q$ is equivalent to minimizing
\[
\Phi_{\mathbf Y}(\mathbf X)=\tfrac12\sum_{k\in\mathcal P}q_k(\mathbf Y)r_k(\mathbf X)^2.
\]
\end{theorem}
\begin{proof}
The function $g_\varepsilon(t)=\sqrt{t+\varepsilon^2}-\varepsilon$ is concave on $[0,\infty)$, with derivative $1/(2\sqrt{t+\varepsilon^2})$. Apply its tangent inequality at $t=r_k(\mathbf Y)^2$ and sum over active pairs.
\end{proof}

\begin{proposition}[Sharp touching scalar majorizer]\label{prop:sharp}
Fix $s\in\R$ and $w>0$. The quadratic
\[
m_s(t)=w\rho_\varepsilon(s)+\frac{w}{2\sqrt{s^2+\varepsilon^2}}(t^2-s^2)
\]
is the unique pointwise smallest quadratic polynomial that majorizes $w\rho_\varepsilon$ on $\R$ and touches it at $s$.
\end{proposition}
\begin{proof}
Let $\pi$ be another such quadratic and put $a=w/\sqrt{s^2+\varepsilon^2}$. Tangency implies $\pi'(s)=w\rho_\varepsilon'(s)=as$. If $s\ne0$, evaluating $\pi(-s)\ge w\rho_\varepsilon(-s)=w\rho_\varepsilon(s)$ gives $\pi''\ge a$. If $s=0$, use $\rho_\varepsilon(t)=t^2/(2\varepsilon)+o(t^2)$ to obtain the same inequality. Therefore
\[
\pi(t)-m_s(t)=\tfrac12(\pi''-a)(t-s)^2\ge0,
\]
and equality of curvature gives equality of the polynomials.
\end{proof}
This is residualwise sharpness in the sense of \citet{de2009sharp}. It does not assert optimality among majorizers that exploit joint constraints on the residual vector. In particular, $Q$ is quadratic in residuals, not in $(\mathbf U,\mathbf B)$.

\section{Conditional SMACOF and Convergence}\label{sec:algorithm}
At outer iteration $v$, set $q_k^{[v]}=q_k(\mathbf X^{[v]})$ and
$\Phi_v(\mathbf X)=\frac12\sum_kq_k^{[v]}r_k(\mathbf X)^2$.

\begin{theorem}[Monotonicity]\label{thm:monotonicity}
If $\Phi_v(\mathbf X^{[v+1]})\le\Phi_v(\mathbf X^{[v]})$ at every iteration, then $\sigma_\varepsilon(\mathbf X^{[v]})$ decreases to a finite nonnegative limit $\bar\sigma$.
\end{theorem}
\begin{proof}
Theorem~\ref{thm:majorization} gives
\[
0\le\sigma_\varepsilon(\mathbf X^{[v+1]})
\le\sigma_\varepsilon(\mathbf X^{[v]})+
\Phi_v(\mathbf X^{[v+1]})-\Phi_v(\mathbf X^{[v]})
\le\sigma_\varepsilon(\mathbf X^{[v]}).
\]
\end{proof}
Monotonicity alone does not imply stationarity. The remaining results use the first exact minimization of the geometric quadratic majorizer.

\subsection{The fixed-weight quadratic update}
For positive active weights $q_{ij}$, let $\mathbf H_q$ be the graph Laplacian with off-diagonal entries $-q_{ij}$ and diagonal entries $\sum_{j\ne i}q_{ij}$. Set
\[
\mathbf G_q=\mathbf V^\top\mathbf H_q\mathbf V,
\qquad \mathbf J=\mathbf1\mathbf1^\top/n.
\]
At $\mathbf Y=(\widetilde{\mathbf U},\widetilde{\mathbf B})$, define the off-diagonal entries of $\mathbf C_q(\mathbf Y)$ by
\begin{equation}
c_{ij}=\begin{cases}
-q_{ij}\delta_{ij}/d_{ij}(\mathbf Y),&d_{ij}(\mathbf Y)>0,\ i\ne j,\\
0,&d_{ij}(\mathbf Y)=0,\ i\ne j,
\end{cases}
\qquad c_{ii}=-\sum_{j\ne i}c_{ij}.
\label{eq:C}
\end{equation}
Inactive pairs have $q_{ij}=0$. The zero-distance convention applies only off the diagonal. The standard Cauchy--Schwarz distance bound yields the touching majorizer
\begin{align}
g_q(\mathbf X\mid\mathbf Y)
={}&\tfrac12\sum_kq_k\delta_k^2
+\tfrac12\tr(\mathbf U^\top\mathbf H_q\mathbf U)
+\tfrac12\tr(\mathbf B^\top\mathbf G_q\mathbf B)\nonumber\\
&-\tr(\mathbf U^\top\mathbf C_q(\mathbf Y)\widetilde{\mathbf U})
-\tr(\mathbf B^\top\mathbf V^\top\mathbf C_q(\mathbf Y)\mathbf V\widetilde{\mathbf B}).
\label{eq:inner_majorizer}
\end{align}
At a zero support distance, the geometric bound is $-\delta_kd_k(\mathbf X)\le0$, so majorization remains valid there.
Under Assumption~\ref{ass:regularity}, the unique minimizer over $\Omega_0$ solves
\begin{align}
(\mathbf H_q+\mathbf J)\mathbf U^+&=\mathbf C_q(\mathbf Y)\widetilde{\mathbf U},\label{eq:Uupdate}\\
\mathbf G_q\mathbf B^+&=\mathbf V^\top\mathbf C_q(\mathbf Y)\mathbf V\widetilde{\mathbf B}.
\label{eq:Bupdate}
\end{align}
Both right-hand sides use the same support point. We solve the positive-definite systems by Cholesky factorization; a truncated pseudoinverse would change the update being analyzed.

\begin{algorithm}[htbp]
\caption{Smoothed-$\ell_1$ conditional SMACOF (sL1-cSMACOF)}\label{alg:sl1}
\begin{algorithmic}[1]
\REQUIRE Active dissimilarities and weights; centered full-rank $\mathbf V$; centered $\mathbf U^{[0]}$ and $\mathbf B^{[0]}$; $\varepsilon>0$; $m\ge1$ inner updates; tolerance $\tau>0$; outer limit $v_{\max}$.
\STATE Set $v=0$ and $s^{[0]}=\sigma_\varepsilon(\mathbf X^{[0]})$.
\WHILE{$v<v_{\max}$}
\STATE Compute $q_k^{[v]}$ from \eqref{eq:weights} and hold them fixed.
\STATE Set $\mathbf X^{[v,0]}=\mathbf X^{[v]}$.
\FOR{$w=0,\ldots,m-1$}
\STATE Form $\mathbf C_{q^{[v]}}(\mathbf X^{[v,w]})$ and solve \eqref{eq:Uupdate}--\eqref{eq:Bupdate} simultaneously to obtain $\mathbf X^{[v,w+1]}$.
\ENDFOR
\STATE Set $\mathbf X^{[v+1]}=\mathbf X^{[v,m]}$ and $s^{[v+1]}=\sigma_\varepsilon(\mathbf X^{[v+1]})$.
\IF{$0\le s^{[v]}-s^{[v+1]}\le\tau\max\{1,s^{[v]}\}$}
\STATE Return $\mathbf X^{[v+1]}$ with objective-tolerance status.
\ENDIF
\STATE Set $v\leftarrow v+1$.
\ENDWHILE
\RETURN The last iterate with iteration-limit status.
\end{algorithmic}
\end{algorithm}
One update ($m=1$) is the default. Repeating exact fixed-weight updates preserves $\Phi_v(\mathbf X^{[v+1]})\le\Phi_v(\mathbf X^{[v,1]})$. More general inner procedures are allowed only when this first-step acceptance condition is satisfied. The asymptotic results refer to the infinite sequence with the finite stopping rule suppressed.

\subsection{Compactness and quantitative first-step control}
\begin{theorem}[Compact sublevels and containment]\label{thm:compact}
Under Assumption~\ref{ass:regularity}, every sublevel of $\sigma_\varepsilon$ in $\Omega_0$ is compact. Global minimizers exist. Any sequence satisfying Theorem~\ref{thm:monotonicity} remains in its initial compact sublevel, has accumulation points, and satisfies $\sigma_\varepsilon(\mathbf X^\ast)=\bar\sigma$ at every accumulation point.
\end{theorem}
\begin{proof}
The set is closed. On a sublevel, coercivity of $\rho_\varepsilon$ bounds all active residuals and hence all active distances. With the base-weight Laplacian $\mathbf H_w$,
\[
\sum_kw_kd_k(\mathbf X)^2
=\tr(\mathbf U^\top\mathbf H_w\mathbf U)+
\tr(\mathbf B^\top\mathbf V^\top\mathbf H_w\mathbf V\mathbf B).
\]
The first quadratic form is positive definite on centered $\mathbf U$, and the second is positive definite by the active-difference span condition. Thus $\mathbf U$ and $\mathbf B$ are bounded. Attainment follows by continuity and compactness. The final assertion follows by continuity along a convergent subsequence.
\end{proof}

Write $\mathbf T_v=\mathbf X^{[v,1]}$ for the first exact update and
\[
E_v=\|\mathbf U^{[v,1]}-\mathbf U^{[v]}\|_F^2+
\|\mathbf B^{[v,1]}-\mathbf B^{[v]}\|_F^2.
\]
\begin{theorem}[First-step decrease and summability]\label{thm:firststep}
Suppose Assumption~\ref{ass:regularity} holds and every accepted outer iterate satisfies
\begin{equation}
\Phi_v(\mathbf X^{[v+1]})\le\Phi_v(\mathbf T_v).
\label{eq:acceptance}
\end{equation}
There is a constant $\mu>0$, independent of $v$, such that
\begin{align}
\sigma_\varepsilon(\mathbf X^{[v]})-\sigma_\varepsilon(\mathbf X^{[v+1]})
&\ge b_v\ge\tfrac\mu2 E_v,\label{eq:decrease}\\
b_v&=\tfrac12\bigl(\|\mathbf U^{[v,1]}-\mathbf U^{[v]}\|_{\mathbf H^{[v]}}^2+
\|\mathbf B^{[v,1]}-\mathbf B^{[v]}\|_{\mathbf G^{[v]}}^2\bigr),\nonumber
\end{align}
where $\|A\|_M^2=\tr(A^\top M A)$. Consequently,
\begin{equation}
\sum_{v=0}^\infty E_v\le\frac{2(\sigma_\varepsilon(\mathbf X^{[0]})-\bar\sigma)}\mu,
\qquad
\min_{0\le v<N}E_v\le\frac{2(\sigma_\varepsilon(\mathbf X^{[0]})-\bar\sigma)}{\mu N}.
\label{eq:stepbound}
\end{equation}
For one inner update, these are bounds on successive outer displacements. For additional inner updates, they concern the first-update displacements.
\end{theorem}
\begin{proof}
Compactness gives $|r_k(\mathbf X^{[v]})|\le R$. Therefore
\[
q_k^{[v]}\ge\frac{w_k}{\sqrt{R^2+\varepsilon^2}}
\ge\theta q_k^{[0]},\qquad
\theta=\frac{\varepsilon}{\sqrt{R^2+\varepsilon^2}}.
\]
We may take
\begin{equation}
\mu=\theta\min\{\lambda_2(\mathbf H^{[0]}),\lambda_{\min}(\mathbf G^{[0]})\}>0,
\label{eq:mu}
\end{equation}
omitting the term involving $\mathbf G^{[0]}$ when $q=0$. Here $g_v=g_{q^{[v]}}(\,\cdot\mid\mathbf X^{[v]})$. Completing squares in \eqref{eq:inner_majorizer} gives $g_v(\mathbf X^{[v]})-g_v(\mathbf T_v)=b_v$. Combining the outer majorization, \eqref{eq:acceptance}, and geometric majorization proves \eqref{eq:decrease}. Telescoping proves \eqref{eq:stepbound}.
\end{proof}
The constant in \eqref{eq:mu} is the smoothed-$\ell_1$ counterpart of the Fair-loss factor $\eta/(R+\eta)$. A computable bound is $R=\sqrt{(S/w_{\min}+\varepsilon)^2-\varepsilon^2}$, where $S=\sigma_\varepsilon(\mathbf X^{[0]})$ and $w_{\min}=\min_kw_k$. No separation assumption has been used in these results.

\subsection{Stationarity and the distance singularity}
\begin{theorem}[Stationarity of accumulation points]\label{thm:stationarity}
Under the assumptions of Theorem~\ref{thm:firststep}, let $\mathbf X^\ast$ be an accumulation point such that
\begin{equation}
d_k(\mathbf X^\ast)>0\quad\text{for every }k\in\mathcal P\text{ with }\delta_k>0.
\label{eq:noncollapse}
\end{equation}
Then $\nabla\sigma_\varepsilon(\mathbf X^\ast)=0$ and $\sigma_\varepsilon(\mathbf X^\ast)=\bar\sigma$. If every accumulation point satisfies \eqref{eq:noncollapse}, every accumulation point is stationary and the gradients of all sufficiently late iterates tend to zero.
At every smooth outer support point,
\begin{equation}
\begin{aligned}
\nabla_{\mathbf U}\sigma_\varepsilon(\mathbf X^{[v]})
 &=\mathbf H^{[v]}(\mathbf U^{[v]}-\mathbf U^{[v,1]}),\\
\nabla_{\mathbf B}\sigma_\varepsilon(\mathbf X^{[v]})
 &=\mathbf G^{[v]}(\mathbf B^{[v]}-\mathbf B^{[v,1]}).
\end{aligned}
\label{eq:gradientidentity}
\end{equation}
If all outer support points are smooth, there is $L<\infty$ such that
\[
\min_{0\le v<N}\|\nabla\sigma_\varepsilon(\mathbf X^{[v]})\|_F^2
\le\frac{2L^2(\sigma_\varepsilon(\mathbf X^{[0]})-\bar\sigma)}{\mu N}.
\]
\end{theorem}
\begin{proof}
At a smooth support point, differentiation gives gradients $(\mathbf H^{[v]}-\mathbf C^{[v,0]})\mathbf U^{[v]}$ and $\mathbf V^\top(\mathbf H^{[v]}-\mathbf C^{[v,0]})\mathbf V\mathbf B^{[v]}$. The first-update normal equations give \eqref{eq:gradientidentity}. Terms with $\delta_k=0$ remain smooth even at $d_k=0$, because their loss is $\sqrt{d_k^2+\varepsilon^2}-\varepsilon$.
Since $q_k^{[v]}\le w_k/\varepsilon$, both matrices in \eqref{eq:gradientidentity} have operator norm at most some $L$, so the gradient norm is at most $L\sqrt{E_v}$.

A subsequence converging to a point satisfying \eqref{eq:noncollapse} eventually lies in a neighborhood where the objective is continuously differentiable. Theorem~\ref{thm:firststep} makes its gradients tend to zero. Gradient continuity proves stationarity. The latent gradient is centered, so ambient stationarity agrees with the constrained condition on $\Omega_0$. The objective value follows from Theorem~\ref{thm:compact}. If all accumulation points are separated, compactness implies an eventual common positive lower bound for the relevant distances; hence the entire tail is smooth and its gradients tend to zero. Squaring the gradient bound and applying \eqref{eq:stepbound} proves the last claim.
\end{proof}

Condition~\eqref{eq:noncollapse} cannot be inferred from smooth residual reweighting or a generic initialization. For $n=2$, $p=q=1$, $\mathbf V=(-1/2,1/2)^\top$, $\delta_{12}=w_{12}=1$, and $\mathbf U=\mathbf B=0$, the zero-distance convention leaves every iterate unchanged. Nevertheless $\mathbf U(t)=(-t/2,t/2)^\top$ decreases the objective for $0<t<1$. The implementation therefore detects positive-dissimilarity collisions and reports a diagnostic stop. Avoiding a finite-iterate collision still does not prove limiting separation.

The software also checks objective descent and records any safeguarded or rejected update. The theorem's exact-update identities are checked only for full exact steps. Small objective change is reported as objective-tolerance termination, not as a prescribed stationarity certificate. Neither the theorem nor these diagnostics establish whole-sequence convergence or global optimality.

\section{Approximation to Exact Least-Absolute Stress}\label{sec:approximation}
Let $F_0(\mathbf X)=\sum_kw_k|r_k(\mathbf X)|$ and $W_\Sigma=\sum_kw_k$.
\begin{proposition}[Uniform envelope and global minimizers]\label{prop:envelope}
For every $\mathbf X$,
\begin{equation}
0\le F_0(\mathbf X)-\sigma_\varepsilon(\mathbf X)\le\varepsilon W_\Sigma.
\label{eq:envelope}
\end{equation}
Under Assumption~\ref{ass:regularity}, both criteria attain their minima. If $\mathbf X_\varepsilon$ is globally optimal for $\sigma_\varepsilon$, then
\[
\min F_0\le F_0(\mathbf X_\varepsilon)\le\min F_0+\varepsilon W_\Sigma.
\]
An additive $\zeta$-optimal smoothed solution has the same upper bound with $\zeta$ added. For any sequence $\varepsilon_j\downarrow0$ and global minimizers $\mathbf X_{\varepsilon_j}$, every accumulation point is globally optimal for $F_0$, and such accumulation points exist.
\end{proposition}
\begin{proof}
The scalar inequalities $|r|\le\sqrt{r^2+\varepsilon^2}\le|r|+\varepsilon$ give \eqref{eq:envelope}. Comparing a smoothed minimizer with an exact-$\ell_1$ minimizer gives the objective bound. The compactness argument of Theorem~\ref{thm:compact} also applies to $F_0$. The minimizers for bounded $\varepsilon_j$ lie in a common compact $F_0$ sublevel, and continuity with \eqref{eq:envelope} proves the final assertion.
\end{proof}
These are global optimization statements, not guarantees for an arbitrary MM output. In particular, the envelope does not convert objective-tolerance stopping into a certificate of near-global least-absolute fitting.

The corresponding scaled Fair loss is $|r|-\eta\log(1+|r|/\eta)$, whose reweighting denominator is $|r|+\eta$. For smoothed $\ell_1$, the denominator is $\sqrt{r^2+\varepsilon^2}$. Both scores are bounded, but equal parameter values do not make their residual weighting identical. Unlike \eqref{eq:envelope}, the Fair-to-absolute difference is unbounded over the full real line at fixed $\eta$. This motivates independent scale calibration in the experiments.

\section{Computational Experiments}\label{sec:experiments}
The study adapts the four application examples, the no-covariate comparison, the tuning analysis, and the iteration study of the Fair-loss rcMDS study \citep{ling2026robust}. All values reported here are computed for the present methods and protocol. We use the same contaminated input and the same start for paired squared-stress, Fair, and smoothed-$\ell_1$ fits within each cell.

\subsection{Protocol and evaluation}
For each replicate, a uniformly selected subset of active undirected pairs is replaced by independent draws from $\mathcal U(2,3)d_{\max}$, where $d_{\max}$ is the largest uncontaminated active dissimilarity. The number of replacements is the nearest integer to the requested fraction of active pairs. The corruption mechanism uses the clean benchmark only to define the experiment. Fitting scales use the \emph{observed} positive active dissimilarities.

Before fitting, divide dissimilarities by their observed positive median $s_\Delta$ and standardize each centered covariate column to unit root-mean-square. In these units, set $\varepsilon=\alpha_{\rm sL1}$ and $\eta=\alpha_{\rm Fair}$; thus their parameters in the original distance units are $\alpha s_\Delta$. Latent starts are centered draws from $\mathcal U(-1/2,1/2)$ and $\mathbf B^{[0]}=0.1\mathbf I$. The maximum is 10,000 outer iterations, with the common relative objective-change rule in Algorithm~\ref{alg:sl1} at $\tau=10^{-9}$. Robust application fits use one inner update. We retain iteration-limit and safeguard outcomes in all summaries.

Let $X$ be a reference configuration and $Y$ the fitted \emph{combined} configuration, padded with zero columns if necessary. Center and Frobenius-normalize both, giving $X_c$ and $Y_c$ with unit norm. The primary shape diagnostic is
\begin{equation}
\mathrm{SS}=1-\|X_c^\top Y_c\|_*^2,
\label{eq:ss}
\end{equation}
which equals symmetric Procrustes error after optimal rotation or reflection and scaling \citep{peres2001well}. Smaller values are better. Mean canonical correlation is reported as a secondary subspace diagnostic. For generated geometries, we additionally use the unaligned distance error
\begin{equation}
\mathrm{NRMSE}_d=
\left\{\frac{\sum_{i<j}(\widehat d_{ij}-d_{ij}^{\rm true})^2}
{\sum_{i<j}(d_{ij}^{\rm true})^2}\right\}^{1/2}.
\label{eq:nrmse}
\end{equation}
All pairs are evaluated, including unobserved pairs in the sparse diagnostic. For perception datasets, reference coordinates are descriptive scales or design variables; their Procrustes scores do not establish recovery of a uniquely identified latent factor.

\subsection{Calibration and comparison methods}
Independent synthetic training replicates select the loss parameters before application evaluation. The smoothed-$\ell_1$ and Fair grids are $\{0.001,0.003,0.01,0.03,0.1,0.3\}$. The Cauchy grid is $\{0.03,0.1,0.3,1\}$. Five training replicates cover the three car-brand conditioning scenarios and contamination levels 0\%, 5\%, and 10\%. For each method, first average SS across scenarios and contamination within replicate. Select the largest grid value whose mean is within one Monte Carlo standard error of the minimum. This is simulation calibration using known training geometry; it is not a data-driven recovery guarantee for a new application.

The unconditional comparison uses squared SMACOF, Fair MDS, smoothed-$\ell_1$ MDS, Cauchy MDS, robust Euclidean embedding (REE), and triangle-filtered MDS (TMDS). The Cauchy criterion is $\frac{c^2}{2}\log(1+(r/c)^2)$, with residual weight $1/(1+(r/c)^2)$, optimized by the same geometric update. REE fits absolute errors in \emph{squared} distances over positive-semidefinite Gram matrices, then truncates the spectrum to the requested dimension \citep{cayton2006robust}. Our projected-subgradient implementation retains the best objective among 2,000 steps with step size $0.1/\sqrt{t}$; it reports a budget outcome, not an SDP optimality certificate. TMDS counts broken triangles, filters edges by the published histogram criterion, and applies weighted SMACOF \citep{blouvshtein2018outlier}. For these small datasets, all triangles are enumerated. If no histogram threshold satisfies both conditions, all edges are retained. A disconnected filtered graph is reported as a failure rather than repaired silently.

These are documented in-house implementations. The Cauchy scale is calibrated in the same conditional training design and then transferred to the unconditional comparisons; it is not separately optimized for those test problems. In the no-covariate study, every method receives the same information and embedding dimension. In the rectangle and face comparisons, the conditional methods receive the stated covariates, while Cauchy, REE, and TMDS do not. Those comparisons therefore combine loss and information effects.

\subsection{Datasets and study sizes}
Synthetic car-brand data use 30 objects and seven attributes: Quality, Safety, Value, Performance, Eco, Design, and Tech. Independent $\mathcal U(0,1)$ attributes are weighted by $(90,88,83,82,81,70,68)/562$. Observed covariates and pairwise dissimilarities receive independent 5\% multiplicative Gaussian noise. Cases 1--3 provide the first four, five, and six attributes, respectively, with the remaining dimension assigned to $\mathbf U$. Ten evaluation replicates are used at 0\%, 2\%, 5\%, 10\%, 15\%, and 20\% replacement contamination. The separate unconditional comparison uses only the first four weighted attributes, a four-dimensional embedding, and ten replicates at 0\%, 5\%, and 10\% contamination.

The kinship data contain 15 terms \citep{rosenberg1975method}; the incomplete Cousin scale is omitted, leaving 14 objects. We condition on Gender or on Gender and Degree, and compare the combined three-dimensional fit with the three recorded scales. Ten paired replicates use 2\% contamination. The rectangle data contain 16 designs \citep{borg1983dimensional}. Width is given and Height is estimated in the first representation; Size (Width plus Height) is given and Shape (Width minus Height) is estimated in the second. Both representations use 100 paired replicates at 5\% contamination. The face data contain 13 expressions and the Pleasantness--Unpleasantness (PU), Attention--Rejection (AR), and Tension--Sleep (TS) scales \citep{abelson1962multidimensional,engen1958dimensional}. Six conditioning subsets, PU, AR, TS, PU--AR, PU--TS, and AR--TS, are evaluated at 0\%, 5\%, and 10\% contamination with ten replicates per cell. The unconditional baselines are compared once per input with the PU-conditioned fit against the same full three-dimensional target.

\subsection{Direct checks of the convergence argument}
The additional generated model uses $n=18$, $p=q=2$, centered Gaussian latent coordinates and covariates, and
$\mathbf B_\star=\left(\begin{smallmatrix}0.8&-0.2\\0.1&0.6\end{smallmatrix}\right)$.
Dissimilarities receive 3\% multiplicative Gaussian noise. We consider the complete graph and a sparse graph formed by independently retaining each pair with probability 0.35 and adding the cycle through all objects. Ten replicates use 0\% or 10\% active-pair contamination. The sL1 method is run with one and five inner updates; squared and Fair fits use one.

For each robust outer iteration, we record $b_v$ in \eqref{eq:decrease}, the actual decrease, the first displacement $\sqrt{E_v}$, the outer displacement, and the residual of \eqref{eq:gradientidentity}. We also record the violation of \eqref{eq:acceptance}, scaled by $\max\{1,2\Phi_v(\mathbf X^{[v]})\}$. Gradient-identity error is divided by $\max\{1,\|\nabla\sigma_\varepsilon\|_F\}$, and decrease violations by $\max\{1,|\sigma_\varepsilon|\}$. Minimum positive-dissimilarity distances, gradient norms, and termination reasons are retained. These are finite-run diagnostics; they do not verify a condition about all limiting configurations.

\FloatBarrier

\subsection{Calibration results and run accounting}\suppressfloats[t]

The calibration rule selected $\alpha_{\rm sL1}=0.1$, $\alpha_{\rm Fair}=0.03$, and $\alpha_{\rm Cauchy}=0.3$. These choices were frozen before the evaluation runs. The complete study contains 720 calibration fits, 2,610 application fits, and 160 theorem-diagnostic fits, for a total of 3,490 attempted fits. Four TMDS inputs yield disconnected filtered graphs and no returned configuration; all other attempts return finite configurations.

\begin{figure}[htbp]\centering
\includegraphics[width=.98\linewidth]{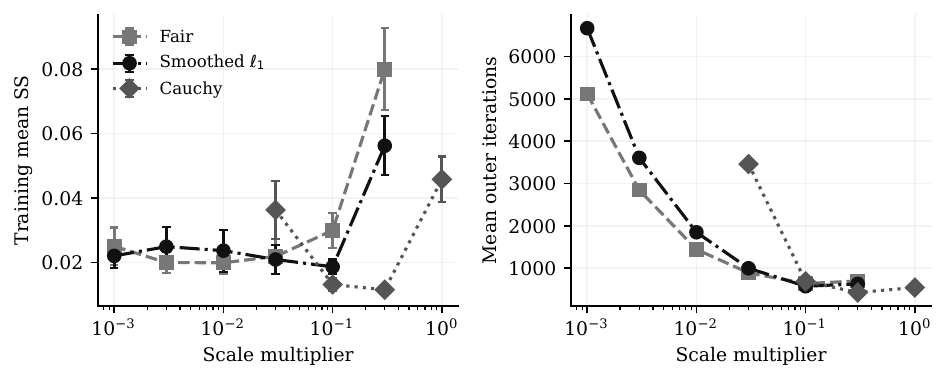}
\caption{Independent training calibration. Left: mean Procrustes SS with one Monte Carlo standard error, computed after within-replicate averaging. Right: mean outer iterations. Parameters are selected separately for each loss.}\label{fig:tuning}
\end{figure}

\begin{table}[htbp]\centering
\caption{Run accounting. Objective tolerance, iteration limits, and the fixed REE budget are distinct termination outcomes.}\label{tab:accounting}
\small
\begin{tabular}{lrrrrr}\toprule
Study & Fits & Tolerance & Limit & REE budget & Failed\\\midrule

Car & 540 & 540 & 0 & 0 & 0\\

Unconditional & 180 & 150 & 0 & 30 & 0\\

Rectangle & 1200 & 997 & 0 & 200 & 3\\

Face & 630 & 586 & 13 & 30 & 1\\

Kinship & 60 & 60 & 0 & 0 & 0\\

Theory & 160 & 160 & 0 & 0 & 0\\

Tuning & 720 & 700 & 20 & 0 & 0\\

\bottomrule\end{tabular}
\end{table}

\FloatBarrier

\subsection{Synthetic recovery results}\suppressfloats[t]

\begin{figure}[htbp]\centering
\includegraphics[width=.98\linewidth]{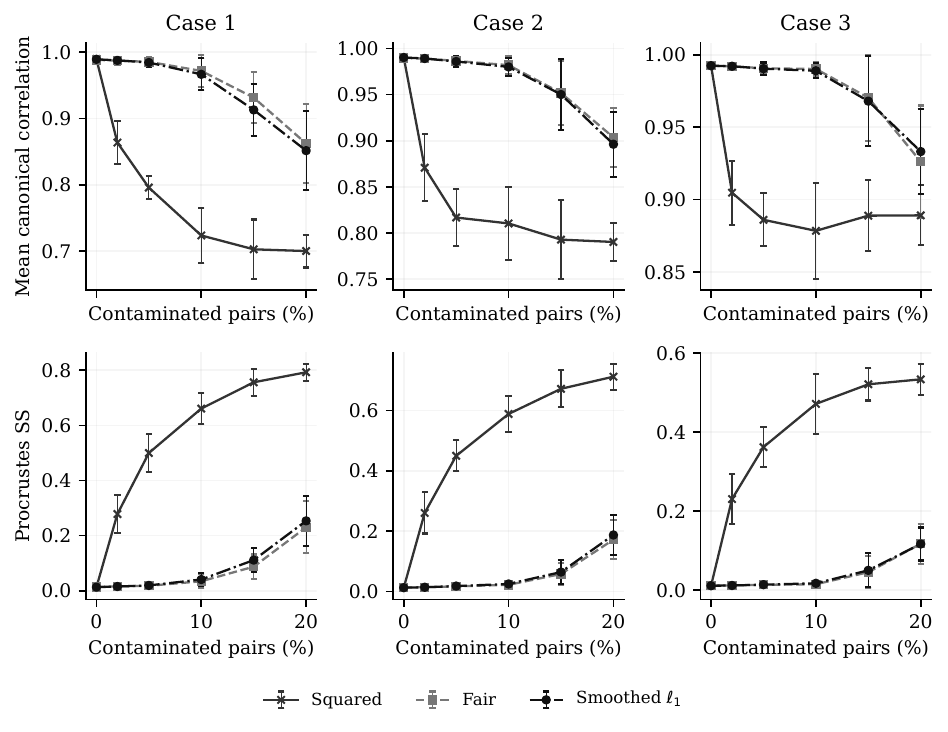}
\caption{Conditional car-brand recovery of the combined seven-dimensional configuration. Columns correspond to four, five, and six given attributes. Top: mean canonical correlation. Bottom: symmetric Procrustes SS. Error bars show Monte Carlo standard deviations across ten replicates.}\label{fig:car}
\end{figure}

Across the three conditional car scenarios, the mean SS values at 0\% contamination are Squared: 0.0122, Fair: 0.0129, sL1: 0.0124. At 20\%, they are Squared: 0.6793, Fair: 0.1732, sL1: 0.1859. Figure~\ref{fig:car} displays each conditioning scenario separately, including variability across paired replicates.

\begin{figure}[htbp]\centering
\includegraphics[width=.98\linewidth]{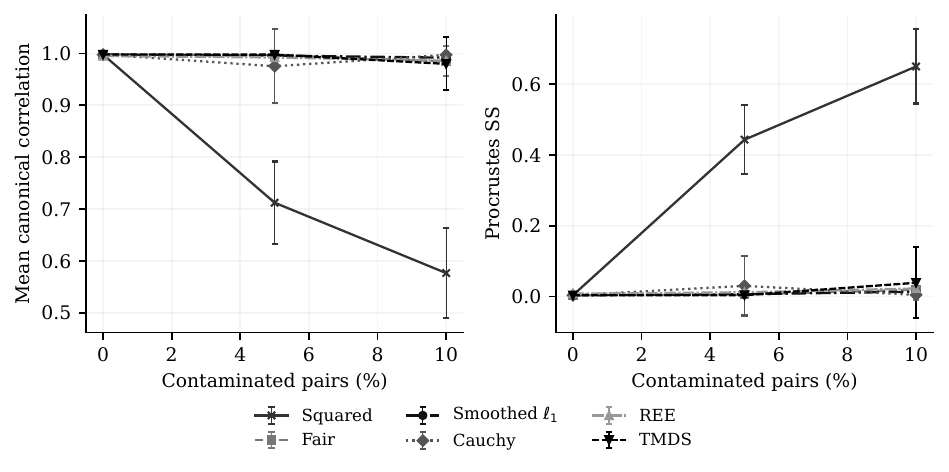}
\caption{Four-attribute car data with no given features. All six methods fit a four-dimensional configuration from the same dissimilarities. Error bars show Monte Carlo standard deviations across ten replicates.}\label{fig:unconditional}
\end{figure}

\begin{table}[htbp]\centering
\caption{Unconditional car recovery at 10\% contamination: means (standard deviations), ten replicates. Distance NRMSE uses original units without alignment.}\label{tab:unconditional}
\small
\begin{tabular}{lrrr}\toprule
Method & SS & Canonical corr. & Distance NRMSE\\\midrule

Squared & 0.6504 (0.1047) & 0.5766 (0.0868) & 0.5931 (0.0341)\\

Fair & 0.0197 (0.0288) & 0.9848 (0.0288) & 0.0495 (0.0188)\\

sL1 & 0.0139 (0.0043) & 0.9917 (0.0027) & 0.0508 (0.0059)\\

Cauchy & 0.0048 (0.0015) & 0.9972 (0.0009) & 0.0271 (0.0023)\\

REE & 0.0229 (0.0058) & 0.9858 (0.0045) & 0.0725 (0.0099)\\

TMDS & 0.0390 (0.1000) & 0.9797 (0.0510) & 0.0491 (0.0599)\\

\bottomrule\end{tabular}
\end{table}

The unconditional comparison holds information access and embedding dimension fixed. REE fits a different squared-distance criterion, and TMDS changes the set of fitted pairs. Their results therefore compare full procedures, not just interchangeable scalar penalties.

At 10\% contamination, Cauchy has the smallest mean SS and distance NRMSE among the six unconditional procedures in Table~\ref{tab:unconditional}. Thus the smoothed loss does not uniformly dominate other robust losses, even when information access is matched.

\FloatBarrier

\subsection{Perception examples}\suppressfloats[t]

\begin{table}[htbp]\centering
\caption{Kinship combined-configuration recovery at 2\% contamination: means (standard deviations) over ten paired replicates. Fourteen complete objects are used.}\label{tab:kinship}
\small
\begin{tabular}{llr}\toprule
Given features & Method & SS\\\midrule

Gender & Squared & 0.7007 (0.0602)\\

Gender & Fair & 0.6981 (0.0859)\\

Gender & sL1 & 0.6606 (0.1063)\\

Gender + Degree & Squared & 0.7225 (0.0483)\\

Gender + Degree & Fair & 0.6327 (0.1707)\\

Gender + Degree & sL1 & 0.6511 (0.1491)\\

\bottomrule\end{tabular}
\end{table}

Kinship SS remains large for every method, despite lower mean error for the robust fits in most conditions. These descriptive scales therefore offer a challenging reference rather than evidence of an accurately recovered latent decomposition.

\begin{table}[htbp]\centering
\caption{Rectangle recovery at 5\% contamination: SS means (standard deviations) over 100 paired replicates for each representation. Conditional methods receive the first design coordinate; Cauchy, REE, and TMDS receive dissimilarities only. TMDS returns 98 width--height and 99 size--shape fits; its means exclude the recorded graph failures.}\label{tab:rectangle}
\small
\begin{tabular}{lrr}\toprule
Method & Width--height & Size--shape\\\midrule

Squared & 0.5101 (0.1344) & 0.5315 (0.1722)\\

Fair & 0.3047 (0.1728) & 0.1971 (0.2131)\\

sL1 & 0.3059 (0.1721) & 0.1893 (0.2137)\\

Cauchy & 0.3346 (0.2231) & 0.1999 (0.1817)\\

REE & 0.2462 (0.0765) & 0.1915 (0.1022)\\

TMDS & 0.3225 (0.2071) & 0.1529 (0.1562)\\

\bottomrule\end{tabular}
\end{table}

Both robust conditional losses reduce mean SS relative to squared stress in the two rectangle representations. Their mean differences from each other are small compared with between-replicate variation. Some unconditional alternatives also attain competitive errors, although they do not use the known design coordinate.

\begin{figure}[htbp]\centering
\includegraphics[width=.98\linewidth]{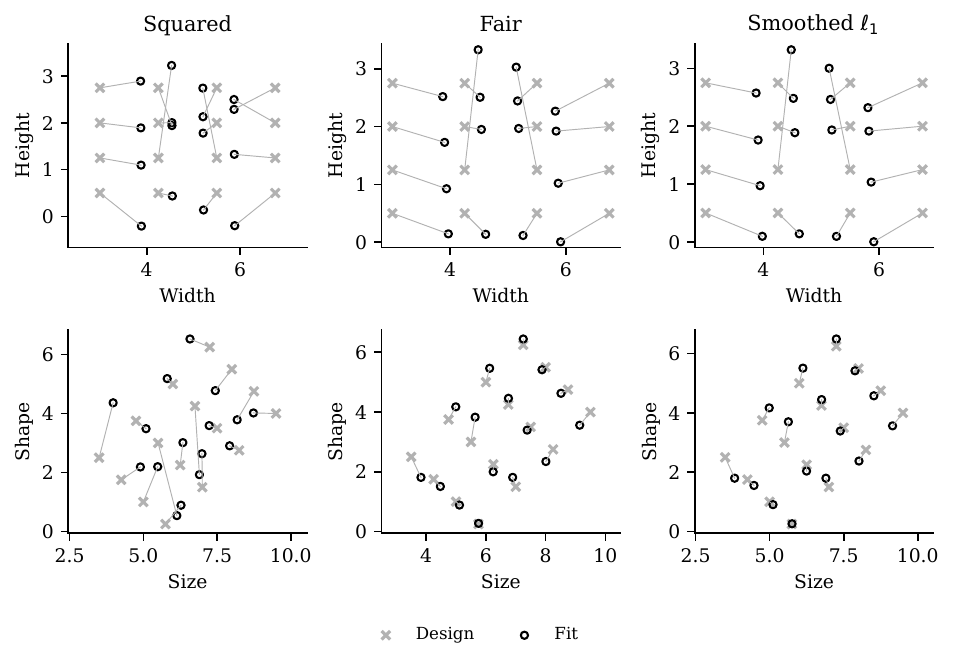}
\caption{Rectangle configurations for the prespecified first replicate at 5\% contamination. Gray crosses are reference designs; open circles are combined fitted configurations after symmetric Procrustes alignment. Lines pair each estimate with its reference.}\label{fig:rectangle}
\end{figure}

\begin{figure}[htbp]\centering
\includegraphics[width=.98\linewidth]{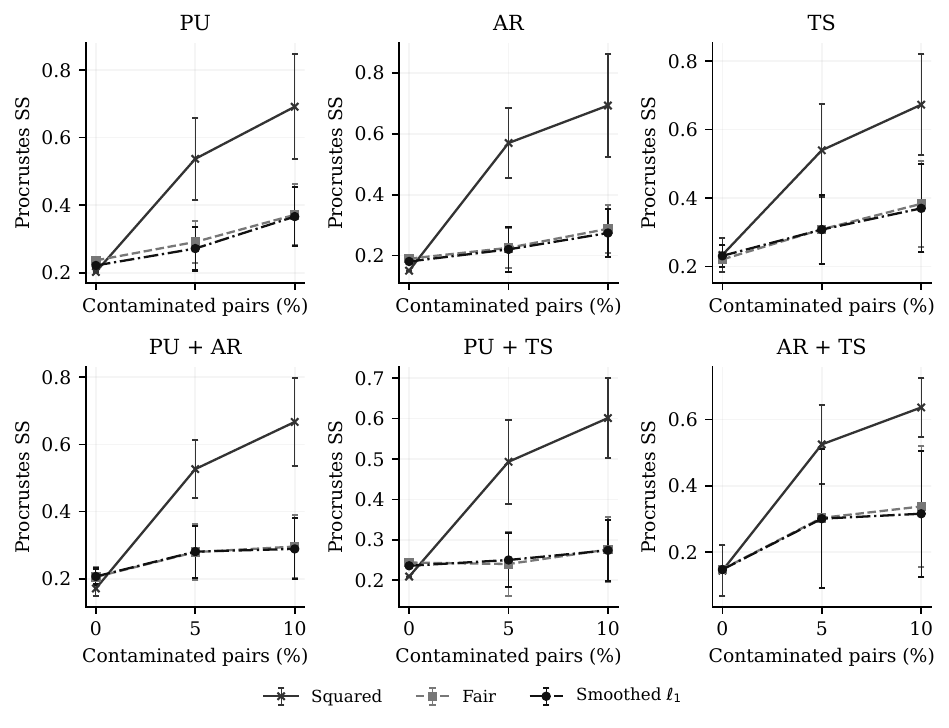}
\caption{Facial-expression recovery against the full PU--AR--TS reference. Panel titles identify given covariates. Curves show SS means and Monte Carlo standard deviations across ten paired replicates.}\label{fig:face}
\end{figure}

For facial expressions, mean SS across the six conditioning patterns at 10\% contamination is Squared: 0.6601, Fair: 0.3253, sL1: 0.3148. Figure~\ref{fig:face} shows the pattern-specific results. The full-target comparison in Table~\ref{tab:facebaseline} uses the PU-conditioned fits and the unconditional baselines on the same inputs.

\begin{table}[htbp]\centering
\caption{Facial-expression full-configuration SS: means (standard deviations). Squared, Fair, and sL1 receive PU; the other methods are unconditional. TMDS has nine returned fits at 10\% because one filtered graph is disconnected; all other cells contain ten fits.}\label{tab:facebaseline}
\small
\begin{tabular}{lrrr}\toprule
Method & 0\% & 5\% & 10\%\\\midrule

Squared & 0.2023 (0.0000) & 0.5372 (0.1220) & 0.6915 (0.1555)\\

Fair & 0.2365 (0.0000) & 0.2912 (0.0623) & 0.3714 (0.0904)\\

sL1 & 0.2215 (0.0000) & 0.2716 (0.0650) & 0.3663 (0.0874)\\

Cauchy & 0.2351 (0.0000) & 0.2613 (0.0637) & 0.2477 (0.0220)\\

REE & 0.2480 (0.0000) & 0.3428 (0.1053) & 0.5919 (0.1661)\\

TMDS & 0.2732 (0.0000) & 0.3889 (0.2095) & 0.4336 (0.1631)\\

\bottomrule\end{tabular}
\end{table}

\FloatBarrier

\subsection{First-step diagnostics and finite-iteration behavior}\suppressfloats[t]

Across the 120 robust diagnostic fits, the largest scaled sufficient-decrease violation is $0$, the largest acceptance violation is $0$, and the largest relative gradient-identity error is $1.96\times10^{-12}$. Every recorded update is a full exact step. The calculation uses the first inner displacement for both one and five updates, as required by Theorem~\ref{thm:firststep}.

\begin{table}[htbp]\centering
\caption{Generated-model distance NRMSE: means (standard deviations) over ten paired replicates. All pairs are evaluated. sL1 columns identify the number of fixed-weight inner updates.}\label{tab:theorydistance}
\small
\begin{tabular}{llrrrr}\toprule
Graph & Contam. & Squared & Fair & sL1 (1) & sL1 (5)\\\midrule

Complete & 0\% & 0.0144 (0.0024) & 0.0302 (0.0508) & 0.0142 (0.0024) & 0.0299 (0.0502)\\

Complete & 10\% & 0.7409 (0.1050) & 0.0629 (0.0702) & 0.0777 (0.0835) & 0.0706 (0.0701)\\

Sparse & 0\% & 0.1161 (0.1306) & 0.0751 (0.0710) & 0.1133 (0.1312) & 0.1133 (0.1312)\\

Sparse & 10\% & 0.9111 (0.1303) & 0.3799 (0.0953) & 0.3771 (0.0905) & 0.3730 (0.0915)\\

\bottomrule\end{tabular}
\end{table}

Five inner updates reduce the reported mean final gradient in each graph and contamination cell, but do not uniformly improve distance recovery. In particular, the clean complete-graph mean NRMSE is larger with five updates than with one. Both variants satisfy the first-step inequalities; the difference illustrates the distinction between descent and recovery in a nonconvex model.

\begin{table}[htbp]\centering
\caption{Final sL1 diagnostics. Gradient norms refer to standardized fitting units. Minimum distances refer to original units and range over positive-dissimilarity active pairs and all ten returned fits. Limits count fits reaching the outer budget.}\label{tab:theorygrad}
\small
\begin{tabular}{llrrrr}\toprule
Graph & Contam. & Inner & Mean gradient & Min. distance & Limits\\\midrule

Complete & 0\% & 1 & 0.0008 & 0.2430 & 0/10\\

Complete & 0\% & 5 & 0.0003 & 0.2430 & 0/10\\

Complete & 10\% & 1 & 0.0049 & 0.2561 & 0/10\\

Complete & 10\% & 5 & 0.0016 & 0.2314 & 0/10\\

Sparse & 0\% & 1 & 0.0005 & 0.2919 & 0/10\\

Sparse & 0\% & 5 & 0.0002 & 0.2919 & 0/10\\

Sparse & 10\% & 1 & 0.0031 & 0.2879 & 0/10\\

Sparse & 10\% & 5 & 0.0014 & 0.2879 & 0/10\\

\bottomrule\end{tabular}
\end{table}

\begin{figure}[htbp]\centering
\includegraphics[width=.98\linewidth]{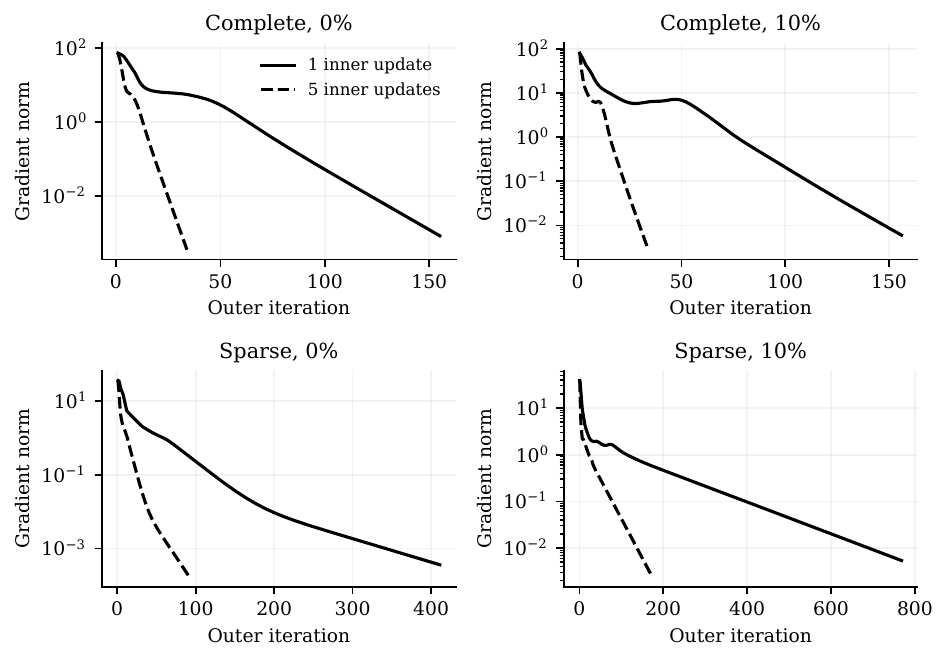}
\caption{Gradient norms at outer support points for the prespecified first replicate. One and five inner updates use the same start and input. Each five-update outer iteration entails more linear-system solves; the horizontal axis is not elapsed computing time.}\label{fig:gradients}
\end{figure}

The exact-$\ell_1$ envelope was also checked on all 80 returned sL1 configurations in this diagnostic. The largest ratio $(F_0-\sigma_\varepsilon)/(\varepsilon W_\Sigma)$ is 0.5260, within the bound in Proposition~\ref{prop:envelope}. This numerical check concerns objective values at returned configurations, not their global optimality. Positive final distances in Table~\ref{tab:theorygrad} are likewise finite-iterate observations, not a proof of limiting separation.

\FloatBarrier

\subsection{Paired contrasts}\suppressfloats[t]

\begin{table}[htbp]\centering
\caption{Paired differences in Procrustes SS. Negative values favor sL1. Intervals are unadjusted 95\% Monte Carlo t intervals. Car, face, and kinship scores are first averaged across conditioning scenarios within each replicate; rectangle contrasts use 100 replicates per representation.}\label{tab:contrasts}
\small
\begin{tabular}{llrr}\toprule
Study & Contrast & Mean & 95\% interval\\\midrule

Car (20\%) & sL1 $-$ Squared & -0.4934 & [-0.5457, -0.4411]\\

Car (20\%) & sL1 $-$ Fair & 0.0127 & [-0.0030, 0.0284]\\

Rectangles: width--height & sL1 $-$ Squared & -0.2042 & [-0.2445, -0.1638]\\

Rectangles: width--height & sL1 $-$ Fair & 0.0012 & [-0.0161, 0.0184]\\

Rectangles: size--shape & sL1 $-$ Squared & -0.3422 & [-0.3894, -0.2949]\\

Rectangles: size--shape & sL1 $-$ Fair & -0.0078 & [-0.0260, 0.0104]\\

Faces (10\%) & sL1 $-$ Squared & -0.3453 & [-0.4126, -0.2780]\\

Faces (10\%) & sL1 $-$ Fair & -0.0105 & [-0.0179, -0.0031]\\

Kinship (2\%) & sL1 $-$ Squared & -0.0557 & [-0.0962, -0.0153]\\

Kinship (2\%) & sL1 $-$ Fair & -0.0095 & [-0.0483, 0.0292]\\

\bottomrule\end{tabular}
\end{table}

All five reported sL1-minus-squared intervals lie below zero. Against Fair loss, only the pooled face contrast excludes zero, with a small advantage for sL1; the other four intervals include zero. These findings support improved recovery over squared stress under the stated contamination designs, without a general preference between the two conditional robust losses.

These paired contrasts quantify uncertainty for the simulated contamination draws and the stated initialization protocol. They do not supply a population-level robustness theorem or a general ordering of loss functions. The termination summaries and gradient diagnostics remain relevant when comparing small differences in recovery.

\section{Discussion}
The analysis separates the roles of residual smoothing, graph connectivity, and distance separation. Smoothing bounds the residual weights above and gives the exact-$\ell_1$ envelope. Connectivity and covariate rank make centered sublevels compact and the quadratic update unique. Limiting separation supplies the differentiability needed for stationarity. None of these statements provides unique identification of the latent and covariate blocks or a guarantee of statistical recovery.

Additional inner updates satisfy the same first-step acceptance condition, but the controlled displacement is $\mathbf X^{[v,1]}-\mathbf X^{[v]}$. It would be incorrect to substitute the full outer displacement into the gradient identity after five updates. More inner work can change both the trajectory and the attained solution; objective decrease alone does not determine which configuration is closer to a reference geometry.

The experiment is limited to small datasets and inflated replacement contamination. The calibration uses synthetic training truth, and its selected values need not transfer to other noise mechanisms. One shared start supports a paired comparison but does not characterize the effect of multistart optimization. REE uses a finite subgradient budget, and TMDS relies on a triangle-based filtering heuristic. Comparisons with unconditional methods in the perception examples also reflect access to covariates. Finally, objective-tolerance and iteration-limit outcomes should be interpreted together with the reported gradients, not as interchangeable stationarity decisions.

\section{Conclusion}
Smoothed least-absolute conditional MDS combines a classical robust penalty with the conditional-distance model and the rcMDS majorization framework. Its residual majorizer is sharp within the class of touching scalar quadratics. A first-exact-update argument establishes compact containment, quantitative descent, and stationarity at separated accumulation points, including accepted outer iterations with additional fixed-weight updates. The uniform loss envelope provides a separate approximation statement for global least-absolute fitting. The computations show improved combined-configuration recovery relative to squared stress in the reported contaminated settings, without a uniform advantage over Fair loss or unconditional robust alternatives. The protocol connects this empirical comparison to conditioning information, independent calibration, and finite-iteration diagnostics. Further work includes collision-aware updates, whole-sequence convergence, broader contamination models, and sparse large-scale solvers.

\section{Disclosure Statement}
Disclosure statement to be supplied by the authors before submission.
\section{Data Availability Statement}
The accompanying archive contains the manuscript and bibliography, base-R solvers, fixed benchmark tables, experiment definitions, random seeds, fitted configurations, raw results, figures, and validation records. All reported values are generated from the included results. A permanent repository identifier and software/data licenses should be supplied before submission.
\begin{samepage}
\bigskip\begin{center}{\large\bfseries SUPPLEMENTARY MATERIAL}\end{center}
\begin{description}
\item[Reproducibility archive:] LaTeX source, R code, benchmark data, saved run checkpoints, numerical summaries, and validation evidence (ZIP).
\item[Source adaptation notes:] A mapping from the supplied Fair-loss manuscript to the smoothed-$\ell_1$ formulas, theorem assumptions, and experiment decisions (plain text).
\end{description}
\end{samepage} 
\bibliography{ref_smoothed_l1}

@article{bui2023dimension,
  title={Dimension reduction with prior information for knowledge discovery},
  author={Bui, Anh Tuan},
  journal={IEEE Transactions on Pattern Analysis and Machine Intelligence},
  volume={46},
  number={5},
  pages={3625--3636},
  year={2024},
  publisher={IEEE}
}

@article{holland1977robust,
  title={Robust regression using iteratively reweighted least-squares},
  author={Holland, Paul W and Welsch, Roy E},
  journal={Communications in Statistics-theory and Methods},
  volume={6},
  number={9},
  pages={813--827},
  year={1977},
  publisher={Taylor \& Francis}
}

@incollection{heiser1988multidimensional,
  title={Multidimensional scaling with least absolute residuals},
  author={Heiser, Willem J},
  booktitle={Classification and Related Methods of Data Analysis},
  pages={455--462},
  year={1988},
  publisher={North-Holland}
}

@article{de2009multidimensional,
  title={Multidimensional scaling using majorization: {SMACOF} in {R}},
  author={De Leeuw, Jan and Mair, Patrick},
  journal={Journal of Statistical Software},
  volume={31},
  pages={1--30},
  year={2009}
}

@article{de2009sharp,
  title={Sharp quadratic majorization in one dimension},
  author={de Leeuw, Jan and Lange, Kenneth},
  journal={Computational Statistics \& Data Analysis},
  volume={53},
  number={7},
  pages={2471--2484},
  year={2009},
  publisher={Elsevier}
}

@article{borg1983dimensional,
  title={Dimensional models for the perception of rectangles},
  author={Borg, Ingwer and Leutner, Detlev},
  journal={Perception \& Psychophysics},
  volume={34},
  number={3},
  pages={257--267},
  year={1983},
  publisher={Springer}
}

@inproceedings{cayton2006robust,
  title={Robust {Euclidean} embedding},
  author={Cayton, Lawrence and Dasgupta, Sanjoy},
  booktitle={Proceedings of the 23rd international conference on machine learning},
  pages={169--176},
  year={2006}
}

@article{peres2001well,
  title={How well do multivariate data sets match? The advantages of a Procrustean superimposition approach over the Mantel test},
  author={Peres-Neto, Pedro R and Jackson, Donald A},
  journal={Oecologia},
  volume={129},
  number={2},
  pages={169--178},
  year={2001},
  publisher={Springer}
}

@article{Bui27052021,
author = {Anh Tuan Bui and Daniel W. Apley},
title = {spc4sts: Statistical process control for stochastic textured surfaces in {R}},
journal = {Journal of Quality Technology},
volume = {53},
number = {3},
pages = {219--242},
year = {2021},
publisher = {Taylor \& Francis},
doi = {10.1080/00224065.2019.1707730},
URL = {https://doi.org/10.1080/00224065.2019.1707730
},
eprint = {https://doi.org/10.1080/00224065.2019.1707730}
}

@article{engen1958dimensional,
  title={The dimensional analysis of a new series of facial expressions.},
  author={Engen, Trygg and Levy, Nissim and Schlosberg, Harold},
  journal={Journal of Experimental Psychology},
  volume={55},
  number={5},
  pages={454--458},
  year={1958},
  publisher={American Psychological Association}
}

@article{abelson1962multidimensional,
  title={Multidimensional scaling of facial expressions.},
  author={Abelson, Robert P and Sermat, Vello},
  journal={Journal of Experimental Psychology},
  volume={63},
  number={6},
  pages={546--554},
  year={1962},
  publisher={American Psychological Association}
}

@article{rosenberg1975method,
  title={The method of sorting as a data-gathering procedure in multivariate research},
  author={Rosenberg, Seymour and Park Kim, Moonja},
  journal={Multivariate Behavioral Research},
  volume={10},
  number={4},
  pages={489--502},
  year={1975},
  publisher={Taylor \& Francis}
}

@article{jaworska2009review,
  title={A review of multidimensional scaling ({MDS}) and its utility in various psychological domains},
  author={Jaworska, Natalia and Chupetlovska-Anastasova, Angelina},
  journal={Tutorials in Quantitative Methods for Psychology},
  volume={5},
  number={1},
  pages={1--10},
  year={2009}
}

@article{charbonnier1997deterministic,
  author  = {Charbonnier, Pierre and Blanc-F{\'e}raud, Laure and Aubert, Gilles and Barlaud, Michel},
  title   = {Deterministic Edge-Preserving Regularization in Computed Imaging},
  journal = {IEEE Transactions on Image Processing},
  volume  = {6},
  number  = {2},
  pages   = {298--311},
  year    = {1997},
  doi     = {10.1109/83.551699}
}

@article{lange2000optimization,
  author  = {Lange, Kenneth and Hunter, David R. and Yang, Ilsoon},
  title   = {Optimization Transfer Using Surrogate Objective Functions},
  journal = {Journal of Computational and Graphical Statistics},
  volume  = {9},
  number  = {1},
  pages   = {1--20},
  year    = {2000}
}

@article{hunter2004tutorial,
  author  = {Hunter, David R. and Lange, Kenneth},
  title   = {A Tutorial on {MM} Algorithms},
  journal = {The American Statistician},
  volume  = {58},
  number  = {1},
  pages   = {30--37},
  year    = {2004},
  doi     = {10.1198/0003130042836}
}

@misc{ling2026robust,
  author={Ling, Xiao and Bui, Anh Tuan},
  title={Robust conditional dimension reduction for dissimilarity data},
  year={2026},
  note={arXiv:2609.06284},
  url={https://arxiv.org/abs/2609.06284}
}

@article{blouvshtein2018outlier,
 author={Blouvshtein, Leonid and Cohen-Or, Daniel},
 title={Outlier Detection for Robust Multi-Dimensional Scaling},
 journal={IEEE Transactions on Pattern Analysis and Machine Intelligence},
 volume={41},number={9},pages={2273--2279},year={2019},
 doi={10.1109/TPAMI.2018.2851513},url={https://arxiv.org/abs/1802.02341}
}
\end{document}